\documentclass[twoside,leqno,10pt, A4]{amsart}
\usepackage{amsfonts}
\usepackage{amsmath}
\usepackage{amscd}
\usepackage{amssymb}
\usepackage{amsthm}
\usepackage{amsrefs}
\usepackage{latexsym}
\usepackage{mathrsfs}
\usepackage{bbm}
\usepackage{amscd}
\usepackage{amssymb}
\usepackage{amsthm}
\usepackage{amsrefs}
\usepackage{latexsym}
\usepackage{mathrsfs}
\usepackage{bbm}
\usepackage{enumerate}
\usepackage{graphicx}
\usepackage{color}
\begin{document}

\newtheorem{theorem}[subsection]{Theorem}
\newtheorem{proposition}[subsection]{Proposition}
\newtheorem{lemma}[subsection]{Lemma}
\newtheorem{corollary}[subsection]{Corollary}
\newtheorem{conjecture}[subsection]{Conjecture}
\newtheorem{prop}[subsection]{Proposition}
\newtheorem{defin}[subsection]{Definition}

\numberwithin{equation}{section}
\newcommand{\mr}{\ensuremath{\mathbb R}}
\newcommand{\mc}{\ensuremath{\mathbb C}}
\newcommand{\dif}{\mathrm{d}}
\newcommand{\intz}{\mathbb{Z}}
\newcommand{\ratq}{\mathbb{Q}}
\newcommand{\natn}{\mathbb{N}}
\newcommand{\comc}{\mathbb{C}}
\newcommand{\rear}{\mathbb{R}}
\newcommand{\prip}{\mathbb{P}}
\newcommand{\uph}{\mathbb{H}}
\newcommand{\fief}{\mathbb{F}}
\newcommand{\majorarc}{\mathfrak{M}}
\newcommand{\minorarc}{\mathfrak{m}}
\newcommand{\sings}{\mathfrak{S}}
\newcommand{\fA}{\ensuremath{\mathfrak A}}
\newcommand{\mn}{\ensuremath{\mathbb N}}
\newcommand{\mq}{\ensuremath{\mathbb Q}}
\newcommand{\half}{\tfrac{1}{2}}
\newcommand{\f}{f\times \chi}
\newcommand{\summ}{\mathop{{\sum}^{\star}}}
\newcommand{\chiq}{\chi \bmod q}
\newcommand{\chidb}{\chi \bmod db}
\newcommand{\chid}{\chi \bmod d}
\newcommand{\sym}{\text{sym}^2}
\newcommand{\hhalf}{\tfrac{1}{2}}
\newcommand{\sumstar}{\sideset{}{^*}\sum}
\newcommand{\sumprime}{\sideset{}{'}\sum}
\newcommand{\sumprimeprime}{\sideset{}{''}\sum}
\newcommand{\sumflat}{\sideset{}{^\flat}\sum}
\newcommand{\shortmod}{\ensuremath{\negthickspace \negthickspace \negthickspace \pmod}}
\newcommand{\V}{V\left(\frac{nm}{q^2}\right)}
\newcommand{\sumi}{\mathop{{\sum}^{\dagger}}}
\newcommand{\mz}{\ensuremath{\mathbb Z}}
\newcommand{\leg}[2]{\left(\frac{#1}{#2}\right)}
\newcommand{\muK}{\mu_{\omega}}
\newcommand{\thalf}{\tfrac12}
\newcommand{\lp}{\left(}
\newcommand{\rp}{\right)}
\newcommand{\Lam}{\Lambda_{[i]}}
\newcommand{\lam}{\lambda}
\newcommand{\af}{\mathfrak{a}}
\newcommand{\sw}{S_{[i]}(X,Y;\Phi,\Psi)}
\newcommand{\lz}{\left(}
\newcommand{\pz}{\right)}
\newcommand{\bfrac}[2]{\lz\frac{#1}{#2}\pz}
\newcommand{\odd}{\mathrm{\ primary}}
\newcommand{\even}{\text{ even}}
\newcommand{\res}{\mathrm{Res}}
\newcommand{\sumn}{\sumstar_{(c,1+i)=1}  w\left( \frac {N(c)}X \right)}
\newcommand{\lab}{\left|}
\newcommand{\rab}{\right|}
\newcommand{\Go}{\Gamma_{o}}
\newcommand{\Ge}{\Gamma_{e}}
\newcommand{\M}{\widehat}
\def\su#1{\sum_{\substack{#1}}}

\theoremstyle{plain}
\newtheorem{conj}{Conjecture}
\newtheorem{remark}[subsection]{Remark}

\newcommand{\pfrac}[2]{\left(\frac{#1}{#2}\right)}
\newcommand{\pmfrac}[2]{\left(\mfrac{#1}{#2}\right)}
\newcommand{\ptfrac}[2]{\left(\tfrac{#1}{#2}\right)}
\newcommand{\pMatrix}[4]{\left(\begin{matrix}#1 & #2 \\ #3 & #4\end{matrix}\right)}
\newcommand{\ppMatrix}[4]{\left(\!\pMatrix{#1}{#2}{#3}{#4}\!\right)}
\renewcommand{\pmatrix}[4]{\left(\begin{smallmatrix}#1 & #2 \\ #3 & #4\end{smallmatrix}\right)}
\def\en{{\mathbf{\,e}}_n}

\newcommand{\ppmod}[1]{\hspace{-0.15cm}\pmod{#1}}
\newcommand{\ccom}[1]{{\color{red}{Chantal: #1}} }
\newcommand{\acom}[1]{{\color{blue}{Alia: #1}} }
\newcommand{\alexcom}[1]{{\color{green}{Alex: #1}} }
\newcommand{\hcom}[1]{{\color{brown}{Hua: #1}} }

\makeatletter
\def\widebreve{\mathpalette\wide@breve}
\def\wide@breve#1#2{\sbox\z@{$#1#2$}%
     \mathop{\vbox{\m@th\ialign{##\crcr
\kern0.08em\brevefill#1{0.8\wd\z@}\crcr\noalign{\nointerlineskip}%
                    $\hss#1#2\hss$\crcr}}}\limits}
\def\brevefill#1#2{$\m@th\sbox\tw@{$#1($}%
  \hss\resizebox{#2}{\wd\tw@}{\rotatebox[origin=c]{90}{\upshape(}}\hss$}
\makeatletter

\title[Bounds for moments of quadratic Dirichlet character sums]{Bounds for moments of quadratic Dirichlet character sums}

%%\date{\today}
\author[P. Gao]{Peng Gao}
\address{School of Mathematical Sciences, Beihang University, Beijing 100191, China}
\email{penggao@buaa.edu.cn}

\author[L. Zhao]{Liangyi Zhao}
\address{School of Mathematics and Statistics, University of New South Wales, Sydney NSW 2052, Australia}
\email{l.zhao@unsw.edu.au}

\begin{abstract}
We establish upper bounds for moments of smoothed quadratic Dirichlet character sums under the generalized Riemann hypothesis, confirming a conjecture of M. Jutila. 
\end{abstract}

\maketitle

\noindent {\bf Mathematics Subject Classification (2010)}: 11L40, 11M06  \newline

\noindent {\bf Keywords}:  quadratic Dirichlet character,  moments

\section{Introduction}\label{sec 1}

  In this paper we are interested in estimating the moments of quadratic Dirichlet character sums given by
\begin{align*}
%%\label{genJacobi}
   S_m(X, Y) :=\sum_{\substack {\chi \in \mathcal S(X) }} \big| \sum_{n \leq Y} \chi(n) \big|^{2m},
\end{align*}  
  where $\mathcal S(X)$ denotes the set of all non-principal quadratic Dirichlet characters of modulus at most $X$. Here $m>0$ is any real number.  \newline

  The case $m=1$ was first studied by M. Jutila \cites{Jutila1} and the best known estimation is given by M. V. Armon \cite[Theorem 2]{Armon}, who showed that 
\begin{align}
\label{S1bound}
   S_1(X, Y) \ll XY(\log X). 
\end{align}  

  For general $m$, a conjecture of Jutila \cite{Jutila2} asserts that for a positive integer $m$, there are constants $c_1(m)$, $c_2(m)$ with values depending on $m$ only,  
\begin{align*}
%%\label{genJacobibound}
   S_m(X, Y) \leq c_1(m)XY^m(\log X)^{c_2(m)}. 
\end{align*}  
  
  In \cites{Virtanen}, H. Virtanen established a weaker version of the above conjecture for the case $m=2$ with the expression $(\log X)^{c_2(m)}$ replaced by
  $X^{\varepsilon}$ for any $\varepsilon>0$.  Other related bounds can be found in \cites{MVa2, Szab}. \newline
  
   It is the aim of this paper to confirm a smoothed version of the above conjecture of Jutila under the assumption of the generalized Riemann hypothesis (GRH). More precisely, we consider a sum of the following form:
\begin{align}
\label{Smoothdef}
   S_m(X,Y; W) := \sum_{\substack{d \leq X \\ (d,2)=1}}\mu^2(d)\Big | \sum_{n}\leg {8d}{n}W \Big(\frac nY \Big )\Big |^{2m},
\end{align}    
   where $W$ is any non-negative, smooth function compactly supported on the set of positive real numbers and $\leg {\cdot}{\cdot}$ denotes the Jacobi symbol. Here we point out (see \cite{sound1}) that the character  $\leg {\cdot}{\cdot}$ is primitive modulo $8d$ for any positive, odd and square-free $d$. Our result gives estimations for $S_m(X,Y; W)$ in terms of the conjectured size.
\begin{theorem}
\label{quadraticmean}
With the notations as above and truth of GRH, we have, for large $X$, $Y$, $\varepsilon>0$ and any real $m \geq 1/2$, 
\begin{align}
\label{mainestimation}
 S_m(X,Y; W) \ll XY^m(\log X)^{m(2m+1)}. 
\end{align}
\end{theorem}

   Using H\"older's inequality and the estimation in \eqref{S1bound}, we can further improve the result in Theorem \ref{quadraticmean} as follows.
\begin{theorem}
\label{quadraticmean1}
With the notations as above and the truth of GRH, we have, for large $X$, $Y$, $\varepsilon>0$ and any real $m \geq 0$, 
\begin{align}
\label{mainestimation1}
 S_m(X,Y; W) \ll
\begin{cases}
  XY^m(\log X)^{m},  \quad 0 \leq m < 1, \\
  XY^m(\log X)^{9m-8},  \quad 1 \leq m < 2, \\
  XY^m(\log X)^{m(2m+1)}, \quad m \geq 2.
\end{cases} 
\end{align}
  Here the above estimation for the case $0 \leq m \leq 1$ holds unconditionally.
\end{theorem}   

  Our proof of Theorem \ref{quadraticmean} is rather simple and makes use of sharp upper bounds on moments of quadratic Dirichlet $L$-functions. Note that this enables our result to be valid for all real $m \geq 1/2$ instead of just positive integers. Our approach here certainly can be applied to bound moments of various other character sums as well. We also point out that upon applying the methods in \cite{Cech2} or \cite{sound1} to evaluate the moment of $|L(s,  \chi_{8d})|^2$ twisted by a Dirichlet character for $\Re(s)=1/2$ and argue as in \cite{Gao2021-2}, one may establish \eqref{mainestimation} for $1/2 \leq m \leq 1$ unconditionally.

\section{Proof of Theorem \ref{quadraticmean}}
\label{sec: mainthm}
 
   We apply the Mellin inversion to obtain that
\begin{align}
\label{charintinitial}
\begin{split}
 S_m(X,Y; W) \ll& 
 \sum_{\substack{d \leq X \\ (d,2)=1}}\mu^2(d)\Big | \int\limits_{(2)}L(s, \chi_{8d})Y^s\widehat W(s) \dif s\Big |^{2m},
\end{split}
\end{align} 	
  where $\widehat W$ stands for the Mellin transform of $W$ given by
\begin{align*}
     \widehat W (s) =\int\limits^{\infty}_0W(t)t^s\frac {\dif t}{t}.
\end{align*}
  Observe that integration by parts implies that for any integer $E \geq 0$,
\begin{align}
\label{whatbound}
 \widehat W (s)  \ll  \frac{1}{(1+|s|)^{E}}.
\end{align}

  Further observe that by \cite[Corollary 5.20]{iwakow} that under GRH, for $\Re(s) \geq 1/2$ and any $\varepsilon>0$,
\begin{align}
\label{Lbound}
 L(s, \chi_{8d}) \ll |ds|^{\varepsilon}. 
\end{align} 
We note here that a bound weaker than \eqref{Lbound} would be sufficient and GRH is not indispensable here. \newline

  The bounds in \eqref{Lbound} and \eqref{whatbound} allow us to shift the line of integration in \eqref{charintinitial} to $\Re(s)=1/2$ to obtain that
\begin{align}
\label{charint}
\begin{split}
 S_m(X,Y; W) \ll & 
 \sum_{\substack{d \leq X \\ (d,2)=1}}\mu^2(d)\Big | \int\limits_{(1/2)}L(s, \chi_{8d})Y^s\widehat W(s) \dif s\Big |^{2m}.
\end{split}
\end{align} 	

 Applying \eqref{whatbound} and H\"older's inequality (note that this requires the condition that $m \geq 1/2$) renders that
\begin{align}
\label{charintreduction}
\begin{split}
 \sum_{\substack{d \leq X \\ (d,2)=1}}\mu^2(d) & \Big | \int\limits_{(1/2)}L(s, \chi_{8d})Y^s\widehat W(s) \dif s\Big |^{2m}\\
  \ll & \sum_{\substack{d \leq X \\ (d,2)=1}}\mu^2(d)\Big (\int\limits_{(1/2)}\Big | \widehat W(s) \Big|^{m/(2m-1)} |\dif s| \Big )^{2m-1}\int\limits_{(1/2)}\Big |L(s, \chi_{8d})Y^s \Big|^{2m} \Big| \widehat W(s)\Big |^{m}|\dif s|  \\
 \ll & Y^m\int\limits_{(1/2)}\sum_{\substack{d \leq X \\ (d,2)=1}}\mu^2(d)\Big | L(s, \chi_{8d})\Big |^{2m} \cdot \Big |\widehat W(s)\Big |^{m} |\dif s|. 
\end{split}
\end{align} 	 

  Now, by \eqref{whatbound} and \eqref{Lbound} again, we may truncate the integral in \eqref{charintreduction} to $|\Im(s)| \leq X^{\varepsilon}$ for any $\varepsilon>0$ with a negligible error. Thus we see from \eqref{charint} and \eqref{charintreduction} that 
\begin{align}
\label{charinttruncation}
\begin{split}
 S_m(X,Y; W) \ll & 
 Y^m\int\limits_{|\Im(s)| \leq X^{\varepsilon}}\sum_{\substack{d \leq X \\ (d,2)=1}}\mu^2(d)\Big | L(s, \chi_{8d})\Big |^{2m} \cdot \Big |\widehat W(s)\Big |^{m} |\dif s|.
\end{split}
\end{align} 	
   
 We modify the proof of \cite[Theorem 2]{Harper} and the proof of \cite[Theorem 2.4]{QShen} (particularly (8.1) of \cite{QShen}) to see that under GRH, we have for $|t| \leq X^{\varepsilon}$ and all $m \geq 0$,  
\begin{align}
\label{upperlower}
 \sum_{\substack{d \leq X \\ (d,2)=1}}\mu^2(d)\Big | L(s, \chi_{8d})\Big |^{2m} \ll_m X(\log X)^{m(2m+1)}.
\end{align}
Note that by applying the argument of Harper \cite{Harper}, one can remove $\varepsilon$-power on the logarithm in the above-mentioned results in \cite{QShen}.    Now upon inserting \eqref{upperlower} into \eqref{charinttruncation}, we immediately obtain the desired result given in \eqref{mainestimation}. This completes the proof of Theorem \ref{quadraticmean}. 

\section{Proof of Theorem \ref{quadraticmean1}}
\label{sec: mainthm1}

   We first note that the estimation given in \eqref{S1bound} is still valid with $S_1(X, Y)$ being replaced by $S_1(X, Y; W)$ for any compactly supported  $W$  by going through the proof of \cite[Theorem 2]{Armon}.  We apply H\"older's inequality to see from this and \eqref{Smoothdef} that for $0 \leq m < 1$, 
\begin{align*}
%%\label{Sest1}
\begin{split}
   S_m(X, Y; W) =& \sum_{\substack{d \leq X \\ (d,2)=1}}\mu^2(d)\Big (1 \cdot  \Big |\sum_{n}\leg {8d}{n}W\Big(\frac nY \Big )\Big |^{2m} \Big )\\
    \leq & \Big (\sum_{\substack{d \leq X \\ (d,2)=1}}\mu^2(d)\Big )^{1-m} \Big (\sum_{\substack{d \leq X \\ (d,2)=1}}\mu^2(d)\Big |\sum_{n}\leg {8d}{n}W \Big(\frac nY \Big )\Big |^{2}\Big )^{m} \ll XY^m(\log X)^{m}.
\end{split}
\end{align*}  
  This gives the estimation for the case $0 \leq m <1$ in \eqref{mainestimation1}. \newline
  
  Similarly, if $1 \leq m < 2$, we deduce from \eqref{S1bound} and \eqref{mainestimation} that for any $p \geq 1$, 
\begin{align*}
%%\label{Sest2}
\begin{split}
   S_m(X, Y; W) =& \sum_{\substack{d \leq X \\ (d,2)=1}}\mu^2(d)\Big (\Big |\sum_{n}\leg {8d}{n}W\Big(\frac nY \Big )\Big |^{2/p} \cdot  \Big |\sum_{n}\leg {8d}{n}W\Big(\frac nY \Big )\Big |^{2m-2/p} \Big )\\
    \leq & \Big (\sum_{\substack{d \leq X \\ (d,2)=1}}\mu^2(d)\Big |\sum_{n}\leg {8d}{n}W\Big(\frac nY \Big )\Big |^{2}\Big )^{1/p}\Big (\sum_{\substack{d \leq X \\ (d,2)=1}}\mu^2(d)\Big |\sum_{n}\leg {8d}{n}W\Big(\frac nY \Big )\Big |^{(2m-2/p)/(1-1/p)}\Big )^{1-1/p}  \\
    \ll & XY^m(\log X)^{1/p+(m-1/p)((2m-2/p)(1-1/p)^{-1}+1)}.
\end{split}
\end{align*} 
  We optimize the exponent of $\log X$ above by setting $1/p=2-m$ to obtain the desired estimation given in \eqref{mainestimation1} for the case $1 \leq m <2$. Note that when $0 \leq m \leq 1$, our estimations above are valid unconditionally since \eqref{S1bound} holds unconditionally. 
 As the case $m \geq 2$ in \eqref{mainestimation1} is just that given in \eqref{mainestimation}, this completes the proof of Theorem \ref{quadraticmean1}.

\vspace*{.5cm}

\noindent{\bf Acknowledgments.}  P. G. is supported in part by NSFC grant 11871082 and L. Z. by the FRG Grant PS71536 at the University of New South Wales.  The authors wish to thank the anonymous referee for his/her careful inspection of the paper and helpful comments.

\bibliography{biblio}
\bibliographystyle{amsxport}

\end{document}